\documentclass[12pt]{article}
\usepackage{graphicx}
\usepackage{amssymb}
\usepackage{amsfonts}
\usepackage{amsthm}
\usepackage{amsfonts}

\newcommand \boundary {\partial}
\newcommand \disc {{\mathbb D}}

\newcommand \no {\noindent}
\newcommand \ra {\rightarrow}

\newcommand{\ba}[1]{\begin{array}{#1}}
\newcommand{\ea}{\end{array}}

\newcommand{\be}{\begin{equation}}
\newcommand{\ee}{\end{equation}}
\newcommand{\bea}{\begin{eqnarray}}
\newcommand{\eea}{\end{eqnarray}}
\newcommand{\beann}{\begin{eqnarray*}}
\newcommand{\eeann}{\end{eqnarray*}}

\def\reff#1{(\ref{#1})}
\newtheorem{conjecture}{Conjecture}

\bibstyle{ams}

\begin{document}

\title{The first order correction to the exit distribution \\
for some random walks}

\author{Tom Kennedy
\\Department of Mathematics
\\University of Arizona
\\Tucson, AZ 85721
\\ email: tgk@math.arizona.edu
}

\maketitle 

\begin{abstract}
We study three different random walk models on several two-dimensional 
lattices by Monte Carlo simulations.
One is the usual nearest neighbor random walk. Another is the 
nearest neighbor random walk which is not allowed to backtrack. 
The final model is the smart kinetic walk.
For all three of these models the distribution of the point where the walk 
exits a simply connected domain $D$ in the plane 
converges weakly to harmonic measure on $\partial D$ 
as the lattice spacing $\delta \ra 0$. 
Let $\omega(0,\cdot;D)$ be harmonic measure for $D$, and  
let $\omega_\delta(0,\cdot;D)$ be the discrete harmonic measure 
for one of the random walk models. Our definition of the random walk 
models is unusual in that we average over the orientation of the lattice
with respect to the domain. 
We are interested in the limit of 
$(\omega_\delta(0,\cdot;D)- \omega(0,\cdot;D))/\delta$.
Our Monte Carlo simulations of the three models lead to the 
conjecture that this limit equals $c_{M,L} \, \rho_D(z)$ times Lebesgue
measure with respect to arc length along the boundary, where  
the function $\rho_D(z)$ depends on the domain, but not
on the model or lattice, and the constant $c_{M,L}$ depends on the model 
and on the lattice, but not on the domain. 
So there is a form of universality for this first order correction.
We also give an explicit formula for the conjectured density $\rho_D$.
\end{abstract}

\newpage

\section{Introduction}

Let $D$ be a simply connected domain in the plane which contains the 
origin. We introduce a lattice with spacing $\delta$ and 
consider a nearest neighbor random walk that starts at the origin. 
We are interested in the point where the walk exits the domain. 
In the limit that the lattice spacing $\delta$ goes to zero,  
the distribution of this random point converges to harmonic
measure which is the distribution of the point where a Brownian motion
exits the domain. In this paper we are interested in the leading order 
term (in $\delta$) in the difference between the random walk exit 
distribution and harmonic measure. 
Our simulations will indicate that this leading order term is 
proportional to $\delta$. So we will refer to it as 
the first order correction. 

In addition to the nearest neighbor random walk, which we will refer 
to as the ordinary random walk, we study two 
other random walks for which the scaling limit of the 
exit distribution is also harmonic measure. 
One model is a nearest neighbor random walk 
that is not allowed to backtrack. In other words, at each step 
the walk picks with equal probability one of its nearest neighbors
other than the site it just came from. 
The third model is the smart kinetic walk, also known as the infinitely
growing self-avoiding walk \cite{weinrib_trugman,kremer1985IGSAW}.
At each time step the walk is only allowed to jump to a nearest neighbor 
that has not been visited before and is not a trapping site
at that time. A site is a trapping site at time $n$ 
if there is no nearest neighbor path from the site to 
the exterior of the domain through sites that are unoccupied 
at time $n$. See \cite{tk_skw}
and references therein for more details.

The convergence of the exit distribution of the ordinary random 
walk to harmonic measure is a well known result. 
An exposition may be found in \cite{lawler_limic}.
A proof of the functional central limit theorem for the random walk without
backtracking on the square lattice may be found in \cite{nbrw}. 
The smart kinetic walk on the hexagonal lattice is equivalent to 
the exploration process for critical percolation on the triangular lattice. 
In this case it has been proved that the exit 
distribution converges to harmonic measure 
\cite{camia_newman, jiang, smirnov,werner2007lectures}. 
For other lattices this is only a conjecture, 
but with good numerical support \cite{tk_skw}.
(The definition of the exploration process in the radial case, i.e., 
between a boundary point and an interior point, is explained in 
section 4.3 of \cite{werner2007lectures}.)

A priori we expect that there are two lattice effects in the first
order correction. One is that  
the lattice introduces a length scale (the lattice spacing).
The other is that the lattice breaks rotation invariance.
Given a domain $D$ we can introduce the lattice with any orientation
with respect to the domain, and the scaling limit will still be 
harmonic measure. However, we expect that the first order 
correction will depend on the orientation of the lattice with respect
to the domain. We attempt to remove this dependence by 
averaging over rotations of the orientation of the lattice 
with respect to the domain. Thus the random walk models that we study are
in fact rotationally invariant, and we expect that the 
main effect of the lattice is the introduction of a length scale. 
There is no doubt that this averaging over the orientation of 
the lattice changes the first order correction in a significant way, 
and so we are not actually studying the first order
correction for the model in which the lattice orientation is kept fixed.
Nonetheless, we see the study of this rotationally averaged model 
as a first step in the study of the first order correction. 

Given a domain $D$ we fix some canonical orientation of the lattice
and then rotate the lattice about the origin
by an angle $\alpha$ with respect to this canonical orientation. 
There is some ambiguity in how one defines the point where the 
random walk exits the domain. The most natural definition
is to linearly interpolate between the steps of the walk so 
that it becomes a piece-wise linear curve in the plane, and we 
can then consider the first point where this curve intersects the 
boundary of the domain. We will refer to the distribution of 
this point as the discrete harmonic measure at $0$ and denote it 
by $\omega^{M,L}_{\delta,\alpha}(0,\cdot;D)$. The symbol
$M$ indicates the model (ordinary
random walk, random walk with no backtracking, or smart kinetic walk),
and $L$ indicates the lattice (square, hexagonal or triangular). 
This is a discrete measure on the finite set of sites where the bonds of the 
lattice with orientation $\alpha$ and spacing $\delta$ intersect the boundary.
Now we average over the orientation by defining 
\bea
\omega^{M,L}_\delta(0,\cdot;D)
= \frac{1}{2 \pi} \int_0^{2 \pi} \omega^{M,L}_{\delta,\alpha}(0,\cdot;D) \, d \alpha
\eea
Because of this averaging over the orientation of the lattice, 
$\omega^{M,L}_\delta(0,\cdot;D)$ 
is a continuous measure on $\partial D$. 
We let $\omega(0,\cdot;D)$ be harmonic measure for $D$. 
We are interested in the difference 
$\omega^{M,L}_\delta(0,\cdot;D)-\omega(0,\cdot;D)$.
We expect it to be of order $O(\delta)$. So we would like to 
compute the limit of 
$[\omega^{M,L}_\delta(0,\cdot;D)-\omega(0,\cdot;D)]/\delta$ as 
$\delta \ra 0$. 

For each of the three models we carry out simulations on 
three lattices - square, hexagonal and triangular. 
Much to our surprise we find that up to an overall constant, the 
first order correction is the same for all three models on all 
three lattices. The constant depends on the model and the lattice, 
but not on the domain. The precise conjecture is as follows.

\begin{conjecture}
For each model $M$ and each lattice $L$ there 
is a constant $c_{M,L}$, and for each simply connected domain
$D$ with smooth boundary there is a function 
$\rho_D(z)$ on the boundary such that 
\bea
\lim_{\delta \rightarrow 0} 
\frac{1}{\delta} [\omega^{M,L}_\delta(0,\cdot;D)
-\omega(0,\cdot;D)] = c_{M,L} \, \rho_D(z) \, |dz|
\eea
Here $|dz|$ is the measure on the boundary given by Lebesgue measure
with respect to arc length. 
The convergence is convergence in distribution. 
The first order difference is a signed measure, 
so the function $\rho_D(z)$ is both positive and negative.
Its integral along the boundary is zero. 

The function $\rho_D(z)$ is determined by the following equation. 
For smooth functions $g(z)$ on the boundary, 
\bea
\int_{\partial D} \, \rho_D(z) \, g(z) \, |dz| =
\int_{\partial D} \, \frac{\partial f}{\partial n}(z) \, \omega(0,|dz|;D)
\eea
where $f$ is the harmonic function in $D$ with boundary values
given by $g$, i.e., it is the solution of Laplace's equation with 
boundary data $g$. The derivative $\frac{\partial f}{\partial n}$
is with respect to the inward unit normal to the domain. 
\end{conjecture}

Again, we should emphasize that the above conjecture is expected to 
be true only if we average over orientations of the lattice. 
If we do not average over the orientations, we expect that the first
order correction depends on the orientation $\alpha$.
We have been deliberately vague about just how smooth the 
boundary of $D$ needs to be in the conjecture since we do not 
have any basis for a precise characterization of the smoothness needed. 

Another possible definition of the random walk exit point from
the domain is to take the first point of the discrete random walk outside the
domain and then orthogonally project this point 
onto the boundary, i.e., take the point on the boundary that is 
closest to the first step the walk takes outside the domain. 
(We will only consider domains with smooth boundary, so when the 
lattice spacing is sufficiently small this point will be unique.) 
With this definition the limiting measure will still
be harmonic measure. However, it is quite possible that 
the first order correction will be different. 
Surprisingly our simulations indicate that it is in fact the same.

In the next section we give a heuristic derivation of the conjecture
for the ordinary random walk. The following section discusses how 
we compute the function $\rho_D(z)$.
Then we discuss the results of our simulations. The paper ends
with some conclusions and open questions. 

\section{Heuristic derivation of the conjecture}

In this section we give a non-rigorous derivation of the conjecture
for the ordinary random walk. As we will see later, this derivation 
only works for the square and triangular lattices. For this derivation 
we use yet another definition of the random walk exit point. 
We run the random walk until it first hits a lattice site in $D$ which 
has a nearest neighbor outside of $D$. 
Then we project this site to the closest site on the boundary of the domain.
The resulting site on the boundary of the domain is the definition of 
the random walk exit point that we will use in this section. 

Let $g$ be a function on the boundary of the simply connected domain $D$.  
Let $f$ be the solution of the continuum Dirichlet problem
\bea
\Delta f(z) =0, \quad z \in  D,  \nonumber \\
f(z)=g(z), \quad z \in \partial D 
\label{dproblem}
\eea
We assume that $\partial D$ and $g$ are sufficiently smooth so that 
$f$ has whatever smoothness is needed below. 
The conjecture says that 
\bea
\lim_{\delta \rightarrow 0} 
\frac{1}{\delta} [\int_{\partial D} g(z) \, \omega_\delta(0,|dz|;D)
- \int_{\partial D} g(z) \, \omega(0,|dz|;D)] 
= c \, \int_{\partial D} \, \frac{\partial f}{\partial n}(z) \, \omega(0,|dz|;D)
\label{hconjecture}
\eea
where $\omega_\delta$ is the exit distribution for one of the 
random walk models and $\omega$ is harmonic measure.

Recall that $\omega_\delta$ is the uniform average over $\alpha$ in 
$[0,2 \pi]$ of $\omega_{\delta,\alpha}$, where $\omega_{\delta,\alpha}$ uses
a lattice at orientation $\alpha$. 
Let $D_{\delta,\alpha}$ be the lattice sites in $D$ for the lattice at orientation 
$\alpha$,  and let $B_{\delta,\alpha}$
be the sites in $D_{\delta,\alpha}$ with at least one nearest neighbor 
outside of $D$. 
For $x \in B_{\delta,\alpha}$ define $g_{\delta,\alpha}(x)$ to be the value of $g$ at 
the point on $\boundary D$ that is closest to $x$. 
(So the line from $x$ to this point is perpendicular to the boundary.) 
Let $f_{\delta,\alpha}$ solve
\bea
\Delta_{\delta,\alpha} f_{\delta,\alpha}(z) 
= 0, \quad z \in D_{\delta,\alpha} \nonumber \\
f_{\delta,\alpha}(z) = g_{\delta,\alpha}(z) \quad z \in B_{\delta,\alpha}
\eea
where $\Delta_{\delta,\alpha}$ is the discrete Laplacian on the lattice 
with spacing $\delta$ at orientation $\alpha$.
The discrete Laplacian is given by 
\bea
\Delta_{\delta,\alpha} h(z) = \frac{a}{\delta^2} \sum_{y:|y-z|=\delta} 
[h(y)-h(z)]
\eea
where $z$ is a lattice site and the sum is over 
the nearest neighbors $y$ of $z$. 
The constant $a$ depends on the lattice and is chosen so that for any 
sufficiently smooth function $h$, 
$\Delta_{\delta,\alpha} h(z)$ converges to $\Delta h(z)$
as $\delta \ra 0$. (On the square lattice $a=1$.)
If we start a nearest neighbor random walk at $0$ and run it 
until it hits $B_{\delta,\alpha}$, then the average of $g_{\delta,\alpha}$ 
with respect to the resulting measure on $B_{\delta,\alpha}$ 
is $f_{\delta,\alpha}(0)$. 
Similarly, $f(0)$ is the integral of $g$ with respect to harmonic measure
on the boundary. 
So \reff{hconjecture} can be rewritten as 
\bea
\lim_{\delta \rightarrow 0} 
\frac{1}{\delta} [\frac{1}{2\pi} \int_0^{2 \pi} f_{\delta,\alpha}(0) \, d\alpha 
- f(0) ] = 
c \, \int_{\partial D} \, \frac{\partial f}{\partial n}(z) \, \omega(0,|dz|;D)
\eea

We now define another solution to the discrete Laplace equation with 
different boundary data. 
Let $\hat{g}_{\delta,\alpha}$ 
be the function on $B_{\delta,\alpha}$ obtained by 
restricting $f$ to $B_{\delta,\alpha}$. 
Let $\hat{f}_{\delta,\alpha}$ solve
\bea
\Delta_{\delta,\alpha} \hat{f}_{\delta,\alpha}(z) = 0, \quad z \in D_{\delta,\alpha}, \nonumber \\
\hat{f}_{\delta,\alpha}(z) = \hat{g}_{\delta,\alpha}(z) \quad z \in B_{\delta,\alpha}
\label{fhatprob}
\eea
Since $f$ is harmonic, if it is sufficiently smooth then  
a simple Taylor series argument shows that $\Delta_{\delta,\alpha} f$ is 
$O(\delta^p)$ where $p=2$ for the square lattice, $p=4$ for the 
triangular lattice, and $p=1$ for the hexagonal lattice.
The discrete maximum principle can then be used to show that 
$f(0)- \hat{f}_{\delta,\alpha}(0)=O(\delta^p)$. (For example, see chapter 1 of 
\cite{numericalpde}.)  
So for the square and triangular lattices which have $p \ge 2$,  
it suffices to study the difference 
$f_{\delta,\alpha}(0) - \hat{f}_{\delta,\alpha}(0)$.
On the hexagonal lattice, the error in replacing 
$f(0)$ by $\hat{f}_{\delta,\alpha}(0)$ is $O(\delta)$ which is as big as the first 
order correction we are trying to derive. So our derivation does not 
work on the hexagonal lattice.
The difference $f_{\delta,\alpha}(0) - \hat{f}_{\delta,\alpha}(0)$
solves the discrete Laplace equation with 
boundary data $g_{\delta,\alpha}-\hat{g}_{\delta,\alpha}$. 
For $z \in B_{\delta,\alpha}$, let $H_{\delta,\alpha}(0,z)$ be the probability 
that the random walk (started at $0$) first hits $B_{\delta,\alpha}$ at $z$. 
Then 
$$
f_{\delta,\alpha}(0) - \hat{f}_{\delta,\alpha}(0)
= \sum_{z \in B_{\delta,\alpha}} \, H_{\delta,\alpha}(0,z) 
[g_{\delta,\alpha}(z)-\hat{g}_{\delta,\alpha}(z)]
$$

Recall that $g_{\delta,\alpha}(z)$ is the value of $g$ at the point on the boundary 
of $D$ which is closest
to $z$. Since $f=g$ on the boundary, we can think of $g_{\delta,\alpha}(z)$ as the 
value of $f$ at this closest boundary point. Recall also that 
$\hat{g}_{\delta,\alpha}(z)$ is the value of $f$ at $z$. 
So a Taylor expansion of the difference of $f$ at these two points yields
\bea
g_{\delta,\alpha}(z)-\hat{g}_{\delta,\alpha}(z) = -\frac{\partial f}{\partial n}(z)
\, dist(z,\partial D) + O(dist(z,\partial D)^2) 
\label{normal_approx}
\eea
Generically, as we move along the boundary the relation of 
the points $z \in B_{\delta,\alpha}$ to the boundary varies. For example, on the 
square lattice some $z$ will have two nearest neighbors that are outside
$D$, while other $z$ will only have one. 
So the probability $H_{\delta,\alpha}(0,z)$ that 
the random walk firsts hits $B_{\delta,\alpha}$ at $z$ will vary significantly 
as we move along $B_{\delta,\alpha}$.

We divide $B_{\delta,\alpha}$ into segments whose length is of 
order $\sqrt{\delta}$. 
As the lattice spacing goes to zero, the number of lattice sites in each 
segment goes to infinity while the length of the segment goes to zero. 
So $\frac{\partial f}{\partial n}$ is essentially constant on each segment. 
Let $\frac{\partial f}{\partial n}(S)$ denote its value on segment $S$.
So we now have 
\bea
f_{\delta,\alpha}(0)-\hat{f}_{\delta,\alpha}(0)  
\approx - \sum_S \frac{\partial f}{\partial n}(S) \,
\sum_{z \in S} \, H_{\delta,\alpha}(0,z) dist(z,\boundary D) \nonumber \\
= - \sum_S \frac{\partial f}{\partial n}(S) \,
\delta(S,\alpha) \, \sum_{z \in S} \, H_{\delta,\alpha}(0,z) \
\eea
where the sum on $S$ is over the segments and 
\bea
\delta(S,\alpha) = \frac{ \sum_{z \in S} H_{\delta,\alpha}(0,z) dist(z,\partial D) }
{\sum_{z \in S} H_{\delta,\alpha}(0,z)}
\label{deltabar}
\eea
Since the boundary is smooth and the lengths of the segments are going to 
zero, the portion of the boundary corresponding to a segment $S$ is
becoming essentially linear. Note that the distance from $z$ to $\partial D$
is of order $\delta$. 
We conjecture that as $\delta$ goes to 
zero, $\delta(S,\alpha)/\delta$ will converge to a limit that depends only 
on the angle of the tangent to the boundary with respect to the lattice 
orientation. Thus when we average over $\alpha$,
$\delta(S,\alpha)$ may be replaced by some constant $c$ times $\delta$. 
For all $\alpha$, ${\sum_{z \in S} H_{\delta,\alpha}(0,z)}$ may be approximated 
by $\omega(0,S;D)$. 
So we now have
\bea
\frac{1}{2\pi} \int_0^{2 \pi} [f_{\delta,\alpha}(0)-\hat{f}_{\delta,\alpha}(0)]
\, d \alpha
\approx - c \delta \sum_S \frac{\partial f}{\partial n}(S) \, \omega(0,S;D)
\eea
After dividing by $\delta$ this converges to 
\bea
-c \int_{\partial D} \, \frac{\partial f}{\partial n}(z) \, \omega(0,|dz|;D)
\eea
So this non-rigorous argument suggests that 
with the definition of the random walk exit point we have used here, 
the constant in the conjecture will be negative.  

We end this section with a probabilistic interpretation of our 
conjecture for the first order correction. Let $A$ be a subarc of the 
boundary. We consider the correction for the probability 
that the random walk exits through $A$, i.e., the probability the random 
walk exits through $A$ minus the probability a Brownian motion exits 
through $A$. Up to an overall constant, our conjecture is 
that to first order in $\delta$ it is given by 
\bea
\delta \int_{\partial D} \frac{\partial f}{\partial n}(z) \, \omega(0,|dz|;D)
\eea
where $f$ is the harmonic function with boundary data $g$ which is 
$1$ on $A$ and $0$ on $A^c$. (By $A^c$ we mean the boundary points that are 
not in $A$.) 
The normal derivative is approximately $(f(z+\delta n)-g(z))/\delta$. 
Since the normal derivative is with respect to the inward normal, 
$z+\delta n$ is a point inside $D$ that is at a distance $\delta$ from 
the boundary. 
Let $\gamma^\delta$ be the curve that is a distance $\delta$ inside $D$ from 
the boundary. In other words, $\gamma^\delta$ is the curve traced by $z+\delta n$
as $z$ runs over the boundary. 
Let $\gamma^\delta_0$ be the part of $\gamma^\delta$ where $g=0$, 
and $\gamma^\delta_1$ the part where $g=1$. Then the above is approximately
\bea
\int_{\gamma^\delta_0} f(z) \, \omega(0,|dz|;D)
+ \int_{\gamma^\delta_1} [f(z)-1] \, \omega(0,|dz|;D)
\eea

Note that 
$f(z)$ is the probability that if we start a Brownian motion
at $z$, then it exits through $A$. 
We denote this by $P(z \ra A)$. 
And $1-f(z)$ is $P(z \ra A^c)$. So the above is 
\bea
\int_{\gamma^\delta_0}   P(z \ra A) \, \omega(0,|dz|;D)
- \int_{\gamma^\delta_1}  P(z \ra A^c) \, \omega(0,|dz|;D)
\eea
The first term is the probability that the Brownian motion gets 
within $\delta$ of the complement of $A$ but then eventually exits 
through $A$. 
The second term is the probability that the Brownian motion gets 
within $\delta$ of $A$  but then eventually 
exits through the complement of $A$.

\section{Computing the conjectured difference}

In \cite{jiang_kennedy} we proved that for a simply connected domain $D$ 
with sufficiently smooth boundary there is a continuous function $\rho_D(z)$ on 
the boundary such that for sufficiently smooth $g$, 
\bea
\int_{\partial D} \, \frac{\partial f}{\partial n}(z) \, \omega(0,|dz|;D)
= \int_{\partial D} \, g(z) \rho_D(z) |dz| 
\eea
where $f$ is the solution of Laplace's equation with boundary data $g$.  
In this section we provide some details on how to explicitly 
compute the function $\rho_D(z)$ so
that we can compare it with the results of our simulations. 

In our simulations we always compute the cumulative distribution
function (CDF) of the exit distribution rather than its density.
So what we need to compute from our conjecture is 
the integral of $\rho_D(z)$ over a subarc of the boundary. 
The five domains that we use in our simulations have the property
that a ray from the origin intersects the boundary in only one point. 
So we can use the polar angle $\theta$ of a point on the boundary to 
parameterize the boundary. 
The following derivation works for any domain with smooth boundary;
one just has to take $\theta$ to be some parameterization of the 
boundary. 
Let $g_\theta(z)$ be the boundary data
that is $1$ on the arc corresponding to 
parameter values in $[0,\theta]$
and is $0$ on the rest of the boundary. Let $f_\theta(z)$ be the 
solution of Laplace's equation with this boundary data. 
We want to compute 
\bea
F_D(\theta)=\int_{\boundary D} \frac{\partial f_\theta}{\partial n}(z) \, 
\omega(0,|dz|;D)
\eea

Let $h(z)$ be the unique conformal map from $D$ to the unit disc which takes 
the origin 
to the origin and takes the point on the boundary of $D$ with polar angle 
$0$ to the point $1$ on the boundary of the unit disc. 
We use $h$ to change the above integral into an integral
over the unit circle. Letting $\phi$ be the angle variable for the unit
circle, $\omega(0,|dz|;D)$ becomes $d \phi/(2 \pi)$ and 
$\frac{\partial f_\theta}{\partial n}(z)$ becomes 
$\frac{\partial f_\theta^\disc}{\partial n}(e^{i \phi}) |h^\prime(h^{-1}(e^{i \phi}))|$.
Here $f_\theta^\disc$ is the harmonic function in the disc with the following
boundary data.  
Let $z(\theta)$ be the point in $\partial D$ with parameter $\theta$, and let 
$\phi(\theta)$ be the polar angle of $h(z(\theta))$. 
The boundary data for $f_\theta^\disc(e^{i \phi})$ is 
$1$ for $\phi=0$ to $\phi=\phi(\theta)$,
and $0$ elsewhere. We now have
\bea
F_D(\theta)= \frac{1}{2 \pi} \int_0^{2 \pi} 
\frac{\partial f_\theta^\disc}{\partial n}(e^{i \phi}) |h^\prime(h^{-1}(e^{i \phi}))| \, 
d \phi
\eea
A straightforward computation using the Poisson kernel for the disc gives
\bea
\frac{\partial f_\theta^\disc}{\partial n}(e^{i \phi})=
\frac{\sin(\phi(\theta)/2)} {2 \pi \sin(\phi/2) \sin((\phi-\phi(\theta))/2)}
\eea
Thus
\bea
F_D(\theta)= \frac{1}{(2 \pi)^2} \int_0^{2 \pi} 
\frac{\sin(\phi(\theta)/2)} {\sin(\phi/2) \sin((\phi-\phi(\theta))/2)}
|h^\prime(h^{-1}(e^{i \phi}))| \, 
d \phi
\eea
We must evaluate this integral numerically. Note that it is only 
conditionally convergent because of the singularities at 
$\phi=0$ and $\phi=\phi(\theta)$.

\section{Simulations}

\begin{figure}[tbh]
\includegraphics{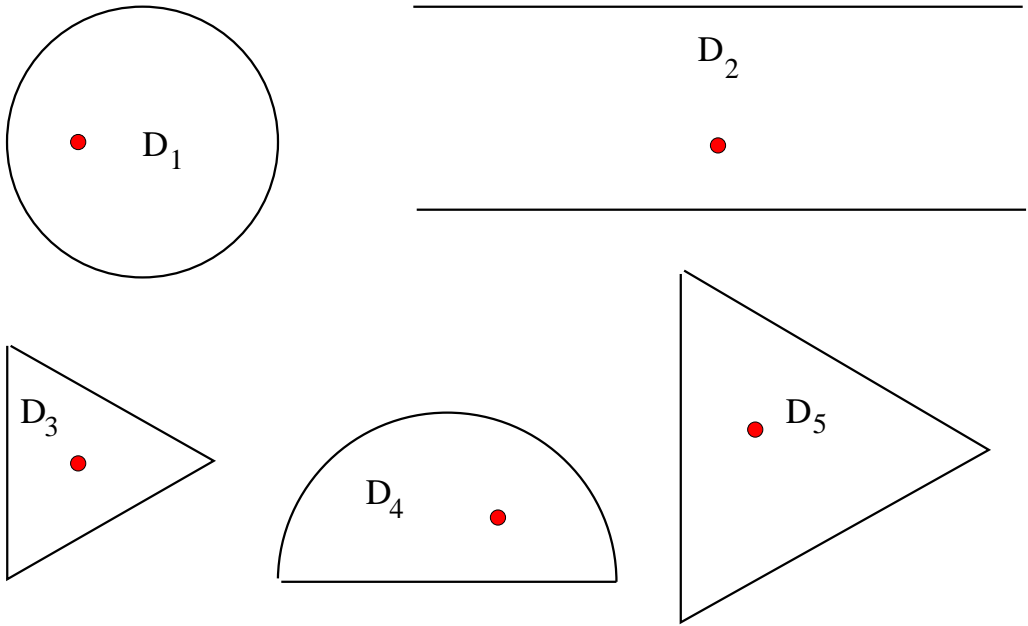}
\caption{\leftskip=25 pt \rightskip= 25 pt 
The five domains. The dot is the origin. 
}
\label{fig_domains}
\end{figure}

We have performed Monte Carlo simulations to compute the exit
distribution for the three different walks (the ordinary random
walk, the random walk with no backtracking and the smart kinetic
walk) on three different lattices (square, hexagonal and triangular) 
in five simply connected domains. 
In all cases the origin is contained in the domain.
The domains are shown in figure \ref{fig_domains} and defined 
as follows.
\beann
D_1 &=& \{ (x,y) : || (x,y)- (1,0)|| < 2 \}, \\
D_2 &=& \{ (x,y) : -1 < y < 2 \}, \\
D_3 &=& T, \\
D_4 &=& \{ (x,y) : || (x,y)- (-3/4,-1)|| < 5/2, \quad y> -1 \}, \\
D_5 &=& \{ (x,y) : (x,y) - (1/2,-1/4) \in T^\prime \}
\eeann
where 
$T$ is the equilateral triangle with vertices
$(2,0), (-1,\sqrt{3})$ and $(-1,-\sqrt{3})$,
and $T^\prime$ is the equilateral triangle with vertices
$(3,0), (-3/2,3 \sqrt{3}/2)$ and $(-3/2,-3 \sqrt{3}/2)$.
We always start the walk at the origin. 
For all three domains the domain has been scaled so that the 
distance from the origin to the boundary is $1$. 

%
%
%
%

We have done simulations with lattice spacings of 
$\delta=0.04, 0.02$ and $0.01$ 
for the three types of walks and the three different lattices
in the five domains. For each case we generated one billion samples.
The average over the orientation $\alpha$ of the lattice is carried out as 
part of the Monte Carlo simulation. To generate each sample we 
first pick an orientation $\alpha$ uniformly from  $[0,2 \pi]$.
Then we run the random walk model on a lattice at that orientation 
until it exits the domain. 

In our simulations we display the results using the cumulative
distribution function (CDF) rather than the density. 
(Computing the density from a simulation requires taking a numerical 
derivative and so adds further uncertainty.)
The differences that we plot are the difference of the CDF for the 
random walk model and the CDF for harmonic measure. 
We plot this difference as a function of the polar angle of the 
boundary point with respect to the origin. 
We will always refer to this function as the difference function 
or simply the difference. Our conjecture for this function is 
$\delta \, c_{M,L} \, F_D(\theta)$ 
where $F_D(\theta)$ is from the previous section.
In the plots we normalize the polar 
angle $\theta$ by dividing it by $2 \pi$. 
For all the simulations described here we define the exit point of 
the random walk by linearly interpolating it between lattice sites and 
looking for the first time this linear interpolation crosses the 
boundary.

\begin{figure}[tbh]
\includegraphics{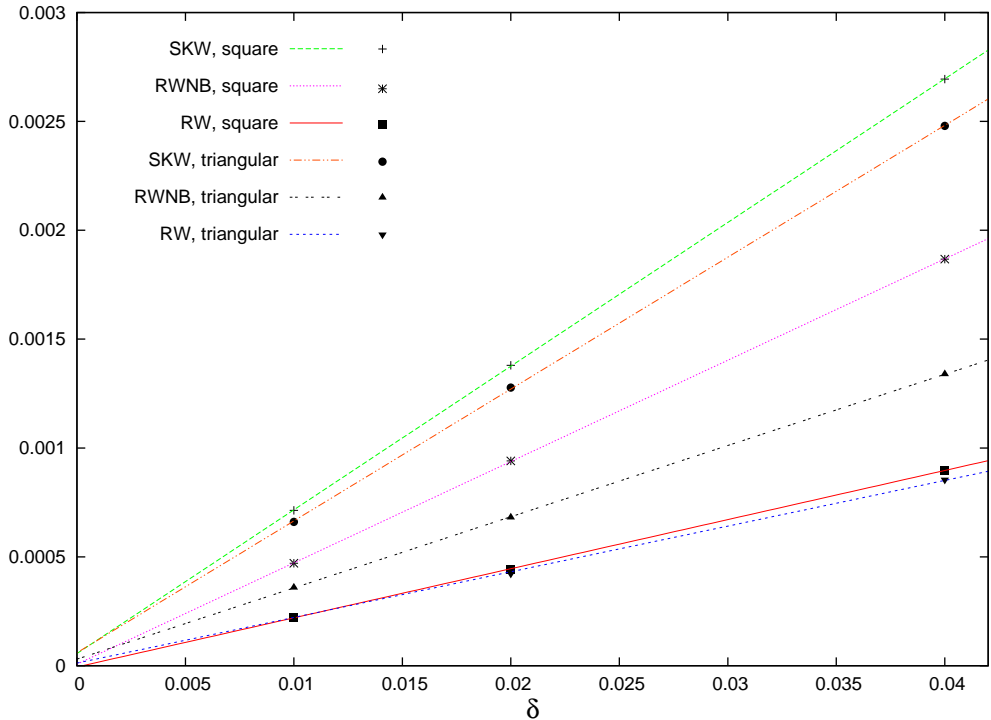}
\caption{\leftskip=25 pt \rightskip= 25 pt 
The $L^1$ norm of the difference function from simulations as a function of the 
lattice spacing $\delta$. The domain is $D_1$.
}
\label{fig_dep_on_delta}
\end{figure}

We first study how the overall magnitude of the difference depends on 
the lattice spacing $\delta$. Our conjecture says it should 
be proportional to $\delta$. In figure \ref{fig_dep_on_delta}
we plot the $L^1$ norm of the difference vs. the lattice spacing
for all three models on the square and triangular 
lattices for domain $D_1$.
The lines shown are fits using linear regression.
The figure supports our conjecture that the magnitude is proportional
to $\delta$. 
The analogous plots for the hexagonal lattice and the other 
domains are similar. 
In the plot RW stands for the ordinary random walk, RWNB stands 
for the random walk without backtracking, and SKW stands for the 
smart kinetic walk.

Next we test the part of the conjecture that asserts the constant
$c_{M,L}$ only depends on the model and the lattice, not on the 
domain. To test this we compute the ratio of 
the $L^1$ norm of the simulation difference function to 
the $L^1$ norm of $\delta \rho_D(z)$.
The conjecture says that this ratio should be independent of the domain. 
Table \ref{table_const} shows this ratio for all three models on 
all three lattices for all five domains. 
If the conjecture is true, then the values in the table should be the 
same in each row. Of course, they are not exactly the same since 
$\delta$ is not quite zero and there are statistical errors from
the simulations. The difference between the smallest and largest
value in a given row ranges from less than $4\%$ to a little more 
than $7\%$. 


\begin{table}
\begin{center}
\begin{tabular}{|l|c|c|c|c|c|c|c|c|c|}
\hline
        &  $D_1$ & $D_2$ & $D_3$ & $D_4$ & $D_5$ \\
\hline
RW hex:    & 0.3996 &  0.3842 &  0.3833 &  0.3885 & 0.3641 \\
RW sq      & 0.3611 &  0.3650 &  0.3704 &  0.3655 & 0.3429 \\ 
RW tri     & 0.3668 &  0.3551 &  0.3626 &  0.3656 & 0.3539 \\ 
RWNB hex   & 1.2108 &  1.1954 &  1.1841 &  1.2051 & 1.1458 \\ 
RWNB sq    & 0.7625 &  0.7825 &  0.7670 &  0.7616 & 0.7376 \\ 
RWNB tri   & 0.5794 &  0.5745 &  0.5651 &  0.5699 & 0.5500 \\ 
SKW hex    & 1.1966 &  1.1710 &  1.1684 &  1.1981 & 1.1400 \\ 
SKW sq     & 1.1632 &  1.1614 &  1.0951 &  1.1316 & 1.0803 \\ 
SKW tri    & 1.0713 &  1.0709 &  1.0150 &  1.0431 & 1.0171 \\ 
\hline
\end{tabular}
\caption{\leftskip=25 pt \rightskip= 25 pt 
Each value in the table is the ratio of 
the $L^1$ norm of the simulation difference function to 
the $L^1$ norm of $\delta \rho_D(z)$.
The row labels indicate the model and the lattice, and the column 
labels indicate the domain. 
}
\label{table_const}
\end{center}
\end{table}

\begin{figure}[tbh]
\includegraphics{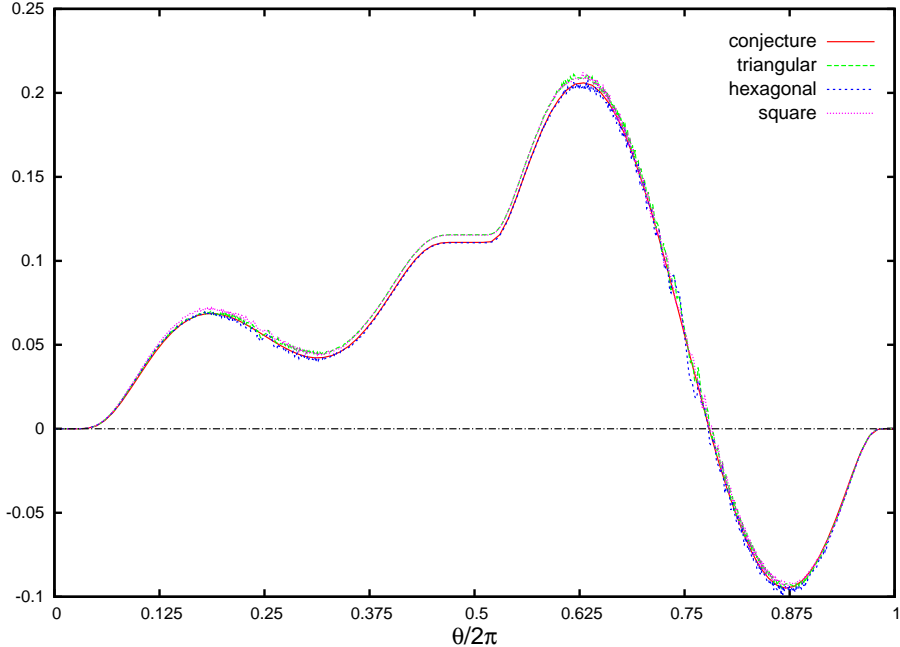}
\caption{\leftskip=25 pt \rightskip= 25 pt 
Rescaled difference for the smart kinetic walk for domain $D_2$.  
}
\label{fig_2_hit_3}
\end{figure}

\begin{figure}[tbh]
\includegraphics{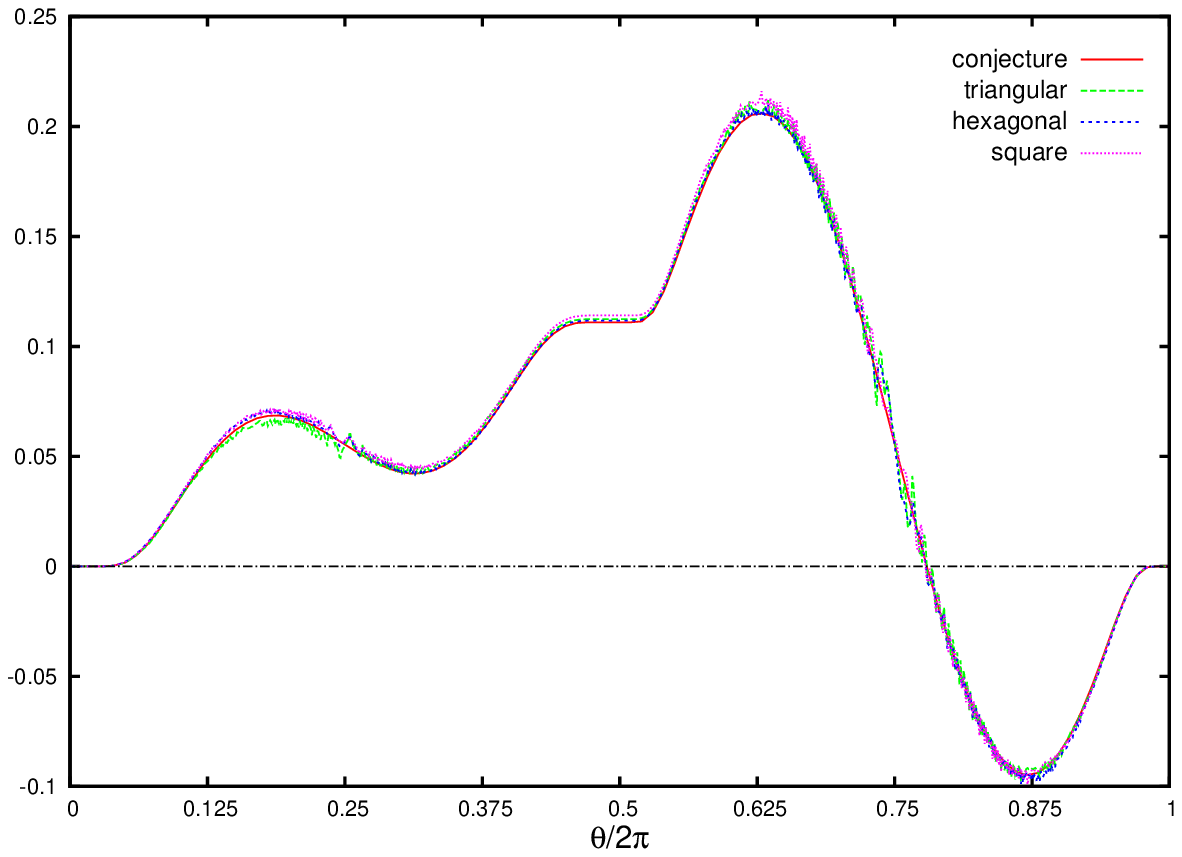}
\caption{\leftskip=25 pt \rightskip= 25 pt 
Rescaled difference for the random walk with no backtracking 
for domain $D_2$.  
}
\label{fig_-2_hit_3}
\end{figure}

\begin{figure}[tbh]
\includegraphics{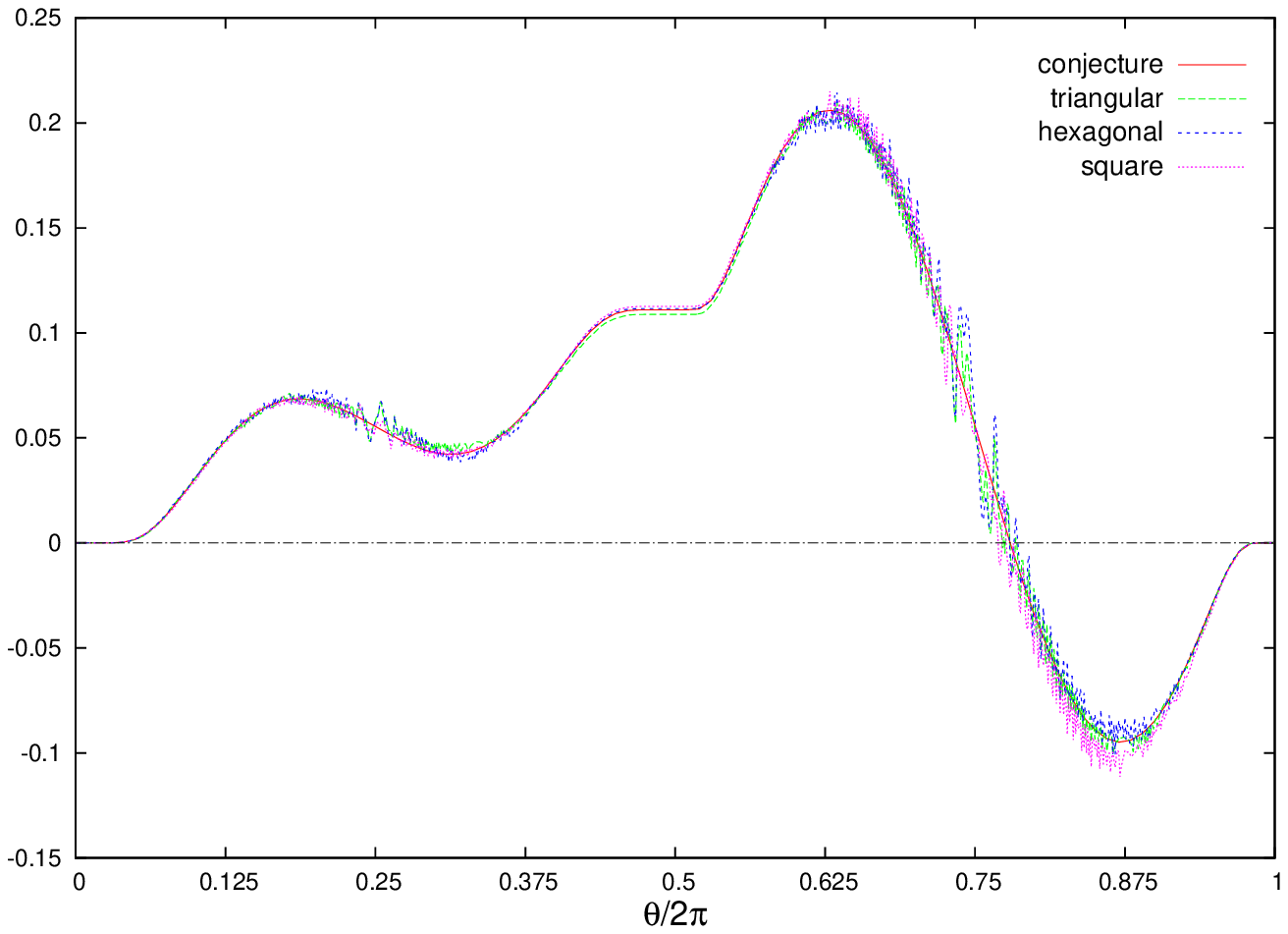}
\caption{\leftskip=25 pt \rightskip= 25 pt 
Rescaled difference for the ordinary random walk for domain $D_2$.  
}
\label{fig_-20_hit_3}
\end{figure}

\begin{figure}[tbh]
\includegraphics{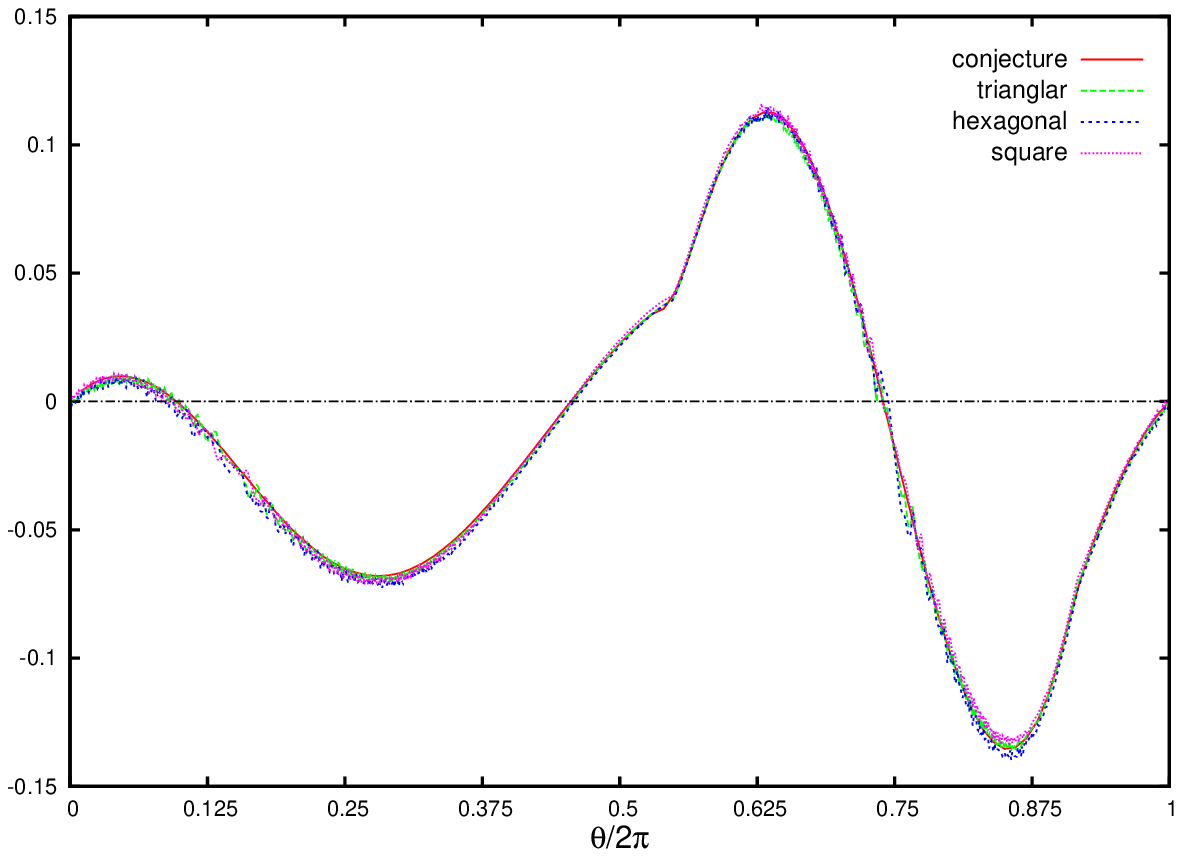}
\caption{\leftskip=25 pt \rightskip= 25 pt 
Rescaled difference for the smart kinetic walk for domain $D_4$.  
}
\label{fig_2_hit_5}
\end{figure}

\begin{figure}[tbh]
\includegraphics{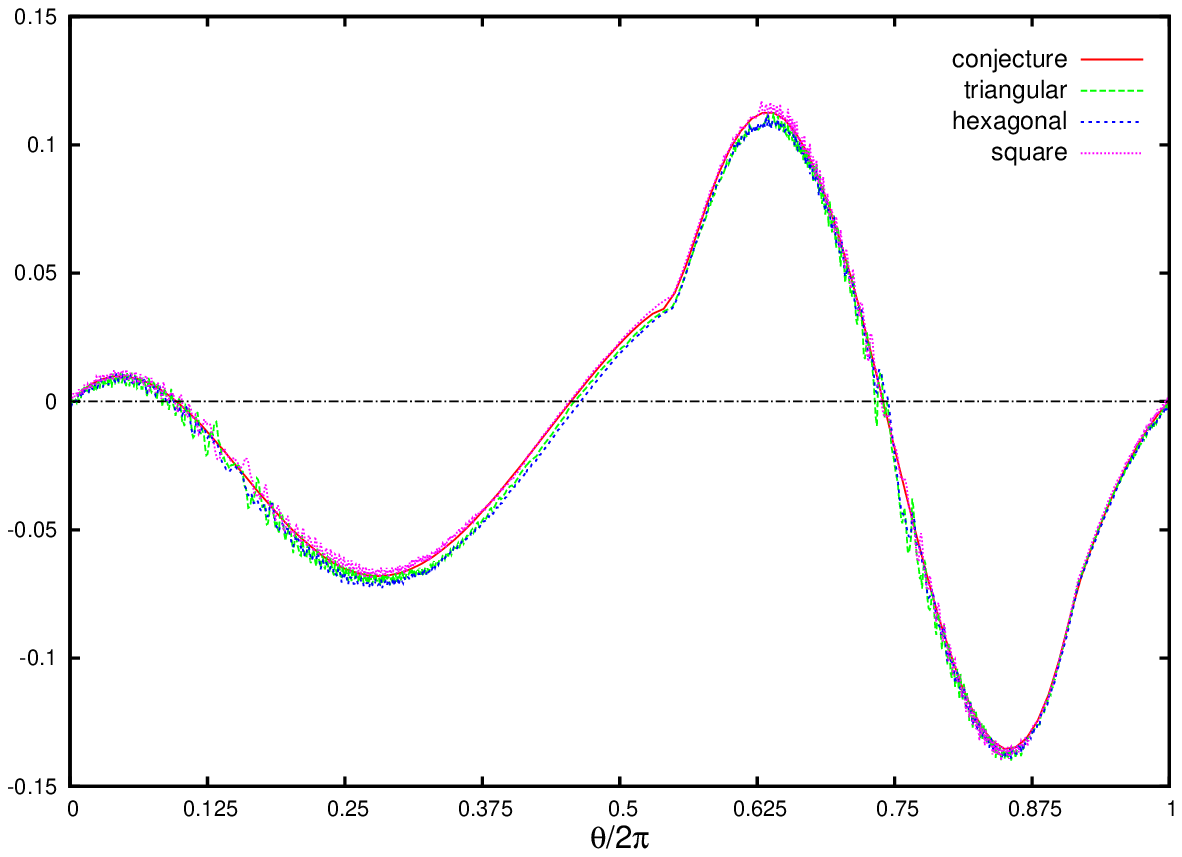}
\caption{\leftskip=25 pt \rightskip= 25 pt 
Rescaled difference for the random walk with no backtracking 
for domain $D_4$.  
}
\label{fig_-2_hit_5}
\end{figure}

\begin{figure}[tbh]
\includegraphics{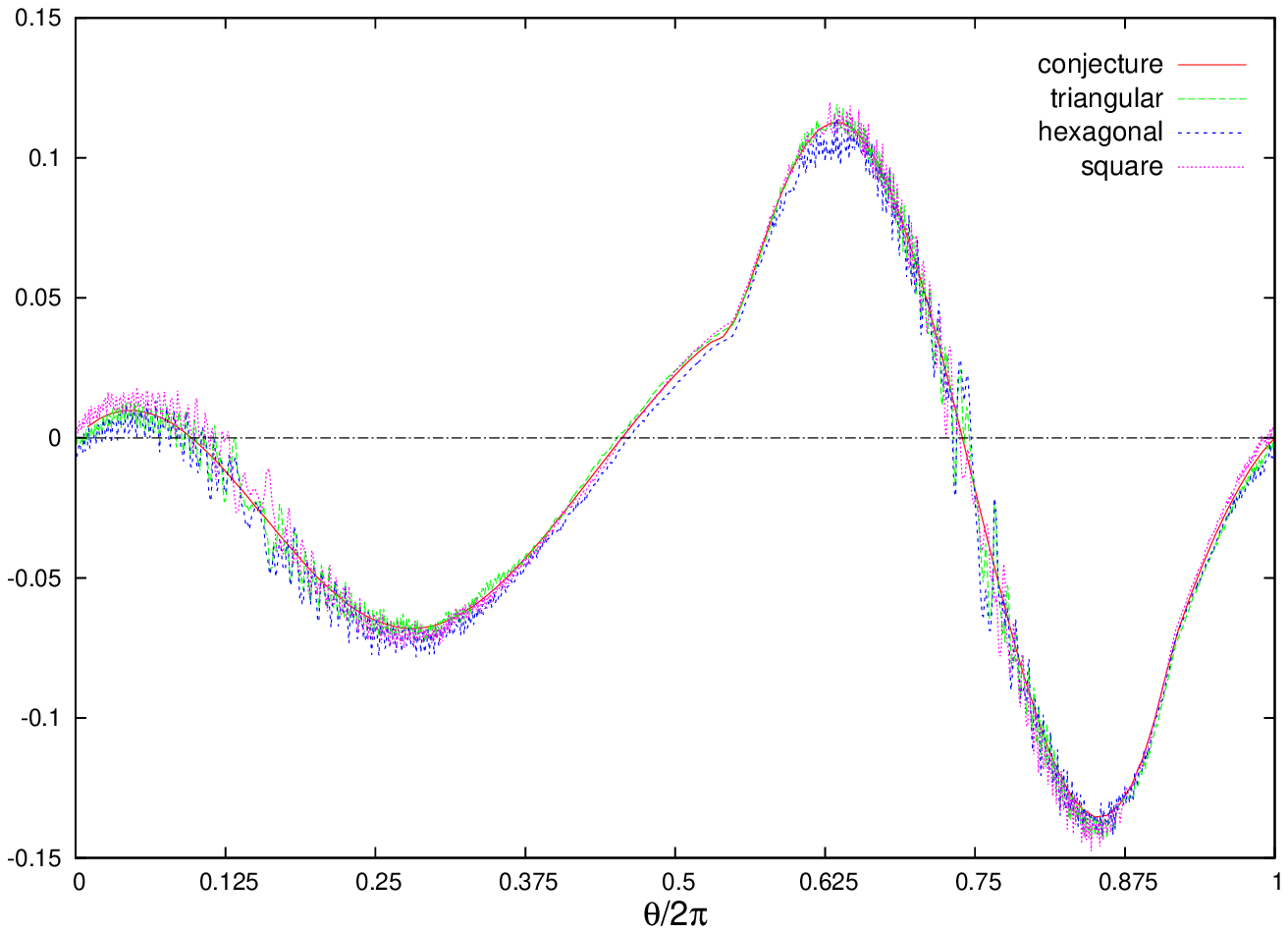}
\caption{\leftskip=25 pt \rightskip= 25 pt 
Rescaled difference for the ordinary random walk for domain $D_4$.  
}
\label{fig_-20_hit_5}
\end{figure}

\begin{figure}[tbh]
\includegraphics{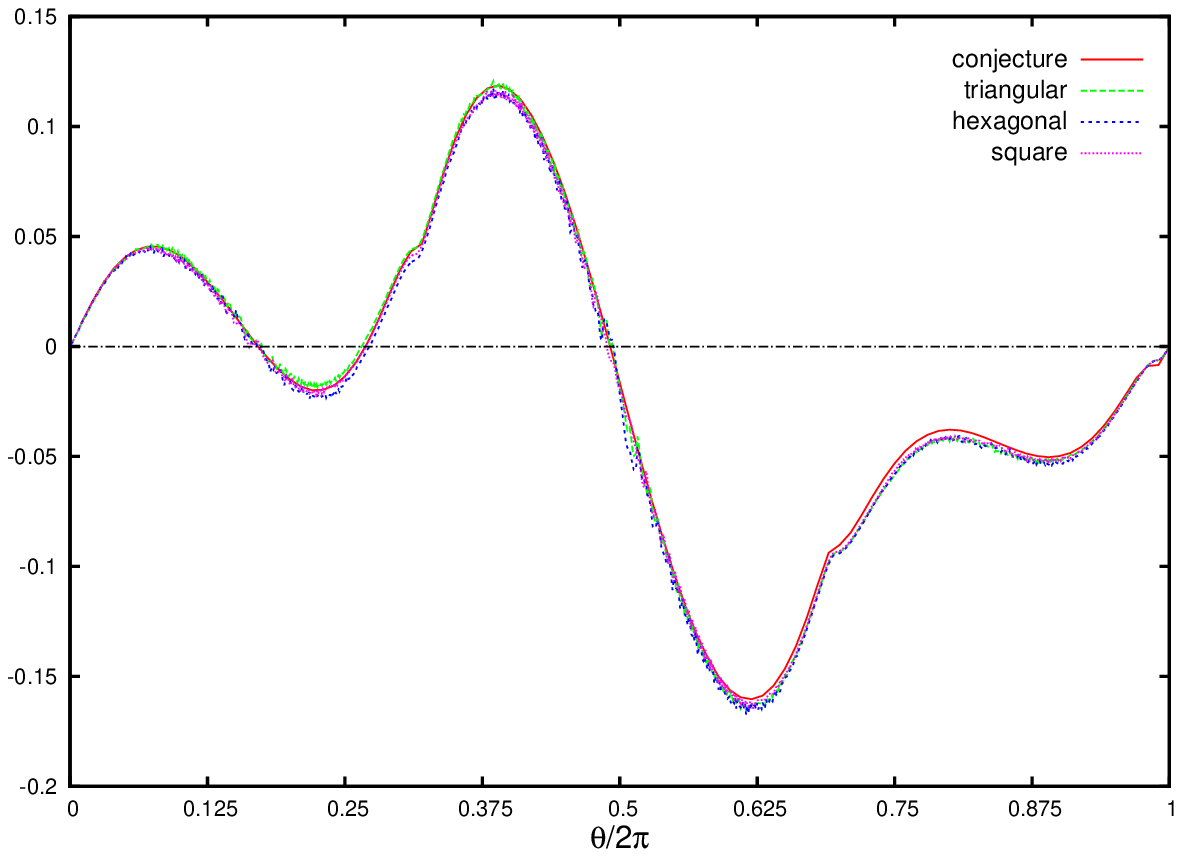}
\caption{\leftskip=25 pt \rightskip= 25 pt 
Rescaled difference for the smart kinetic walk for domain $D_5$.  
}
\label{fig_2_hit_6}
\end{figure}

\begin{figure}[tbh]
\includegraphics{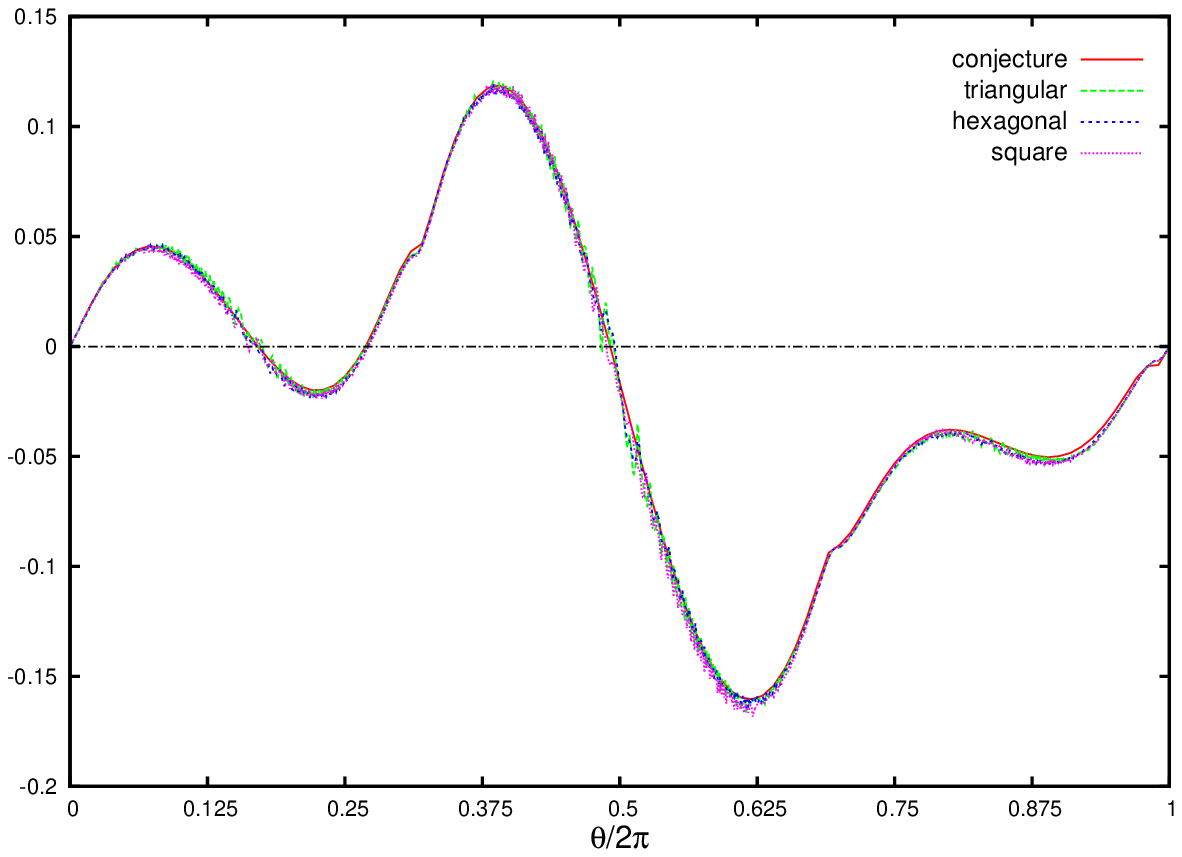}
\caption{\leftskip=25 pt \rightskip= 25 pt 
Rescaled difference for the random walk with no backtracking 
for domain $D_5$.  
}
\label{fig_-2_hit_6}
\end{figure}

\begin{figure}[tbh]
\includegraphics{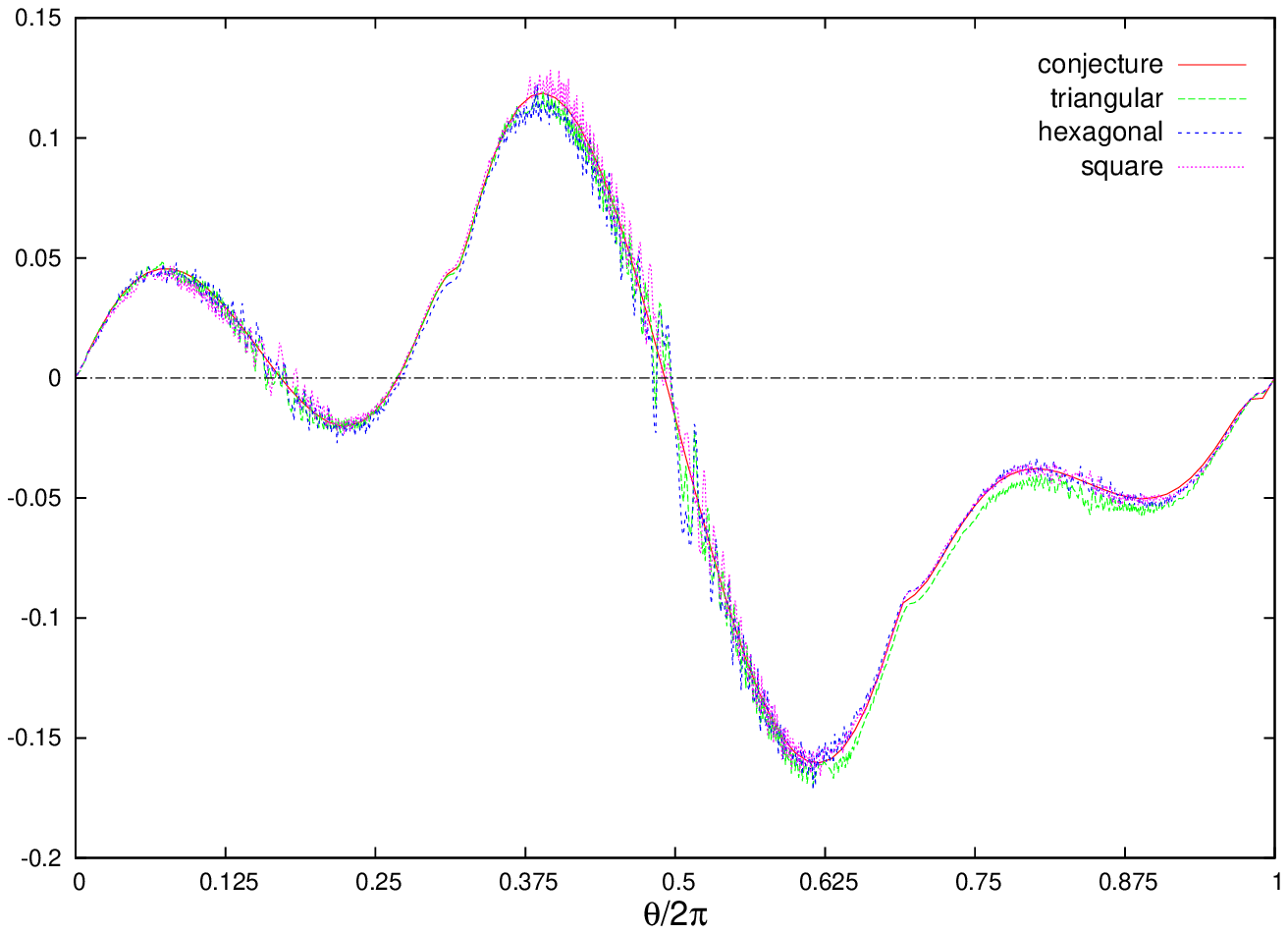}
\caption{\leftskip=25 pt \rightskip= 25 pt 
Rescaled difference for the ordinary random walk for domain $D_5$.  
}
\label{fig_-20_hit_6}
\end{figure}

Finally we compare the polar angle dependence of the difference functions from
our simulations with the function predicted by our 
conjecture. To do this we rescale all the difference functions 
by dividing them by $\delta \, c_{M,L}$. The conjecture says that 
for a given domain $D$, all nine of these rescaled difference 
functions should be equal to the function $F_D(\theta)$.
For the values of $c_{M,L}$ we 
use the average of the values in the table over the five domains. 
We only consider the simulations with the smallest lattice spacing
($\delta=0.01$) for this comparison. 
Since there are three models being simulated on three lattices, 
for each domain we have nine difference functions in addition to the 
function from the conjecture. Rather than attempt to show all ten of 
these functions in a single plot, we plot the three difference functions
for a single model on all three lattices in a single plot. 
For reasons of space we only show the results for the three 
domains $D_2, D_4$ and $D_5$.
The agreement between the simulations and the conjecture for the other two 
domains is at least as good as the agreement for these three. 
The resulting nine figures are figures \ref{fig_2_hit_3} to 
\ref{fig_-20_hit_6}. Each figure contains four curves, but it 
can be difficult to distinguish them since they are so close together. 

\begin{figure}[tbh]
\includegraphics{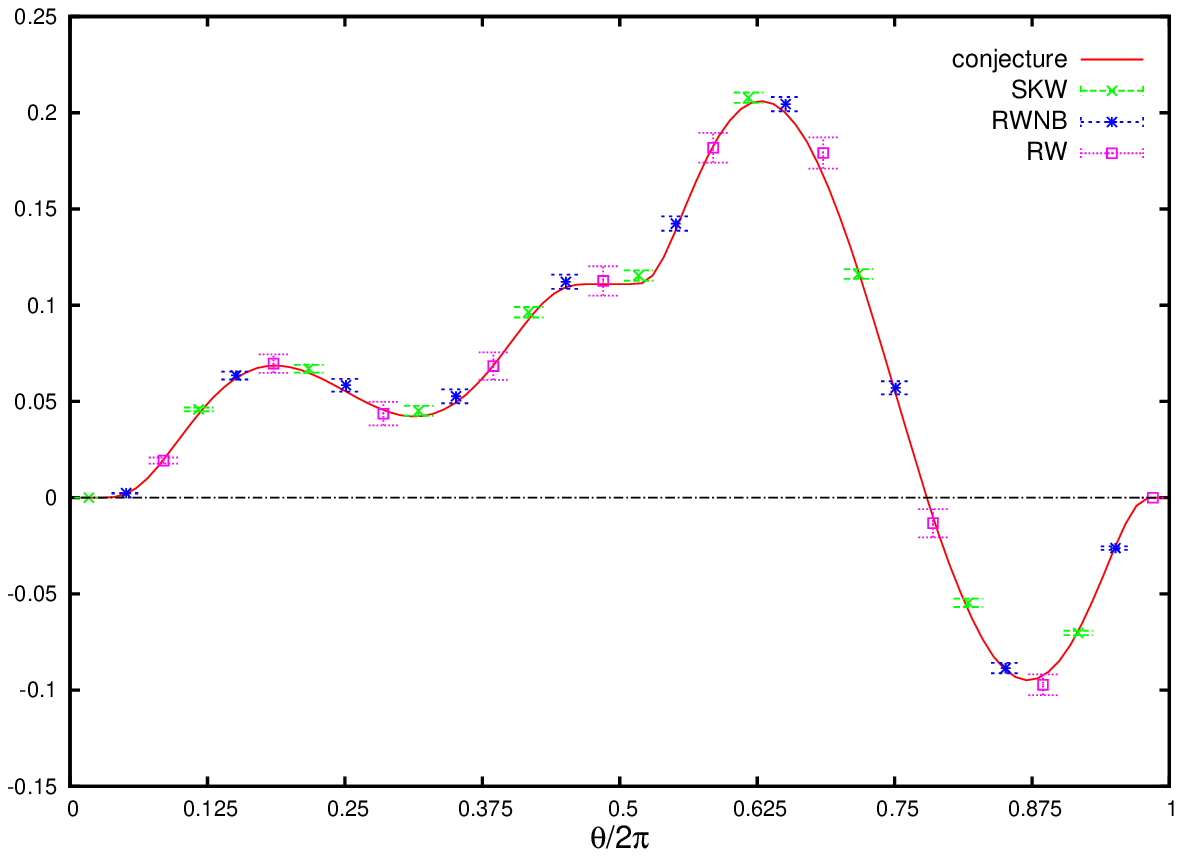}
\caption{\leftskip=25 pt \rightskip= 25 pt 
Error bars for the three models on the square lattice. The domain
is $D_2$.
}
\label{fig_error_bars}
\end{figure}

We have not shown any error bars in figures \ref{fig_2_hit_3} to 
\ref{fig_-20_hit_6} to keep these figures simple. 
In figure \ref{fig_error_bars} we show a few error bars for the 
simulations of the three models on the square lattice for domain 
$D_2$. The curve shown is the conjectured difference function.
The error bars for the other lattices and domains are similar
in size. In this figure the error bars have been rescaled 
by dividing them by $\delta \, c_{M,L}$ just as in the previous 
figures. Because of this rescaling the error bars for the smart kinetic 
walk appear smaller than those of the random walk with no backtracking,
which in turn appear smaller than those of the ordinary random walk.
Before the rescaling the error bars are comparable in size. 
These error bars only represent the error from the Monte Carlo simulation,
i.e., from the fact we can only generate a finite number of samples. 
There are also errors from the fact that the lattice spacing is 
nonzero. They are not included in the error bars in the figure.

\section{Conclusions}

We have studied numerically the exit distribution of 
the ordinary nearest neighbor random walk, the nearest neighbor 
random walk with no backtracking and the smart kinetic walk.
We modify the usual definition of these models by also averaging
over the orientation of the lattice with respect to the domain. 
For all three models the limiting distribution as the lattice 
spacing goes to zero is harmonic measure.
Our simulations support the conjecture that for all three models 
on the three lattices studied, the difference between the random 
walk exit distribution and harmonic measure is, to first order 
in the lattice spacing $\delta$, given by $\delta \, c_{M,L} \, \rho_D(z) |dz|$
where the constant $c_{M,L}$ depends only on the lattice and the model, and 
the density function $\rho_D(z)$ depends only on the domain $D$. 
Thus there is a sort of universality for this first order 
correction. 

The most obvious open problem is to prove this conjecture; the 
ordinary random walk model is the natural model to consider first. 
In \cite{jiang_kennedy} the authors considered a random walk which takes 
place in the continuum rather than on a lattice. 
The steps of the walk are i.i.d. random variables which are 
uniformly distributed over a disc of radius $\delta$.
Note that this model is rotationally invariant. 
It was proved that the conjecture holds for this model and an 
explicit expression for the constant which is the analog of $c_{M,L}$
was given. This proof depended heavily on detailed estimates of the 
Green's function for the model near the boundary of the domain.
For the ordinary random walk the needed estimates on the Green's function
appear to be much harder.

Our non-rigorous derivation of the conjecture for the ordinary random 
walk was based on the generator of this walk being the lattice
Laplacian. The random walk with no backtracking is also Markovian
if we enlarge the state space to include the direction the walk just 
came from. This suggests that it may be possible to 
give a non-rigorous derivation of the conjecture for this model 
using the generator.
By contrast the smart kinetic walk is far from being Markovian, 
so a non-rigorous derivation for the conjecture for this model will
require a new approach.

In principle the non-rigorous derivation of the conjecture for 
the ordinary random walk gives a way to compute the constant, 
eq.  \reff{deltabar}. It should be possible to compute this 
expression by Monte Carlo simulation and compare the result
with the corresponding values in the table.

Perhaps the most important question for future research is to determine
the first order correction when we do not average over the orientation
of the lattice. If one fixes the orientation of the lattice, then 
the exit distribution is a discrete distribution, and this discreteness
complicates the Monte Carlo computation of 
the first order correction in the convergence
to harmonic measure. To avoid this one can average the orientation 
of the lattice only over a subinterval of $[0,2 \pi]$. 
Preliminary simulations show that the difference function
when we average over a subinterval is significantly different 
from the difference function when we average over the full interval. 
We would like to generalize our conjecture to this case of 
averaging over a subinterval. In particular, it would be very 
interesting to determine if the first order correction still 
has a universal shape in the sense that the density function is 
the same for all three models on all three lattices except for an 
overall constant of proportionality.

\begin{appendix}

\section{Error in earlier versions}

The first two versions of this paper contained
an error in the simulations of the ordinary random walk on the 
hexagonal lattice. The error has been corrected in this version. 
There are two changes in the results after this correction. 
The most significant change is in the values in the first row of 
table \ref{table_const}. The corrected values are roughly the 
same size as the values in the second and third rows of the table 
for the random walk on the square and triangular lattices. 
The incorrect values for the hexagonal lattice were much larger. 
The other change is in the differencess for the ordinary random walk 
on the hexagonal lattice which are plotted in figures 
\ref{fig_-20_hit_3}, \ref{fig_-20_hit_5},\ref{fig_-20_hit_6}. 
The overall magnitudes of these differences change
significantly with the correction,
but once these differences are rescaled, it is virtually impossible 
to distinguish the correct and incorrect plots. 

The error in the simulation for the ordinary random walk on the hexagonal
lattice was as follows. At each step the walk should have three equally probable
choices : backtrack to the site it just came from, turn left by 
$60$ degrees or right by $60$ degrees.  The error was a sign error 
so that when the walk should have backtracked it instead took another step 
in the direction it had just taken. 
It is rather surprising that even with the error the rescaled differences
exhibit the same universal behavior seen in all the other models. 
Our understanding of this is as follows. With the error, the simulation is 
actually simulating a kind of random walk on the triangular lattice. 
At each step the walk picks with equal probability one of 
three possible steps:
go forward, turn left by $60$ degrees, turn right by $60$ degrees. 
In the scaling limit the exit distribution of this model will be harmonic 
measure. The simulations indicate that our conjecture holds for this model
as well, with a different constant $c_{M,L}$ from the other models. 

\end{appendix}

\bigskip
\bigskip

\no {\it Acknowledgments:} 
This research was partially supported by NSF grant DMS-1500850.
An allocation of computer time from the UA Research Computing High Performance 
Computing (HPC) and High Throughput Computing (HTC) 
at the University of Arizona is gratefully acknowledged.
The author thanks Jianping Jiang for many 
stimulating conversations about this research.


\bigskip

\begin{thebibliography}{}

\bibitem{camia_newman}
F. Camia, C. M. Newman, 
Critical percolation exploration path and SLE$_6$ : a proof of convergence.
Probab. Theory Related Fields
{\bf 139},473--519 (2007).
Archived as arXiv:math/0605035 [math.PR].

\bibitem{nbrw}
R. Fitzner, R. van der Hofstad,
Non-backtracking random walk, 
J. Stat. Phys. {\bf 150}, 264-284 (2013).
Archived as arXiv:1212.6390 [math.PR].

\bibitem{jiang} 
J. Jiang,
Exploration processes and SLE$_6$. Preprint (2014). 
Archived as arXiv:1409.6834 [math.PR].

\bibitem{jiang_kennedy} J. Jiang, T. Kennedy, 
The difference between a discrete and continuous harmonic measure,
J. Theoret. Probab., to appear. 
Archived as arXiv:1506.04313 [math.PR].

\bibitem{tk_skw} T. Kennedy, 
The Smart Kinetic Self-Avoiding Walk and Schramm-Loewner Evolution, 
J. Stat. Phys. {\bf 160}, 302-320 (2015).
Archived as arXiv:1408.6714 [math.PR].

\bibitem{kremer1985IGSAW}
K. Kremer, J. W. Lyklema, 
Indefinitely growing self-avoiding walk.
Phys. Rev. Lett. {\bf 54}, 267 (1985).

\bibitem{lawler_limic}
G. F. Lawler, V. Limic,
{\it Random Walk: A Modern Introduction}.
Cambridge University Press (2010).

\bibitem{smirnov}
S. Smirnov, Critical percolation in the plane: Conformal invariance, 
Cardy's formula, scaling limits.
C. R. Math. Acad. Sci. Paris
{\bf 333}, 239--244 (2001).
Archived as arXiv:0909.4499 [math.PR].

\bibitem{weinrib_trugman}
A. Weinrib, S. A. Trugman, 
A new kinetic walk and percolation perimeters.
Phys. Rev. B  {\bf 31}, 2993 (1985).

\bibitem{werner2007lectures}
W. Werner,
Lectures on two-dimensional critical percolation, 
{\it Statistical Mechanics (IAS/Park City mathematics series v. 16)}, 
S. Sheffield, T. Spencer (eds.) (2007).
Archived as arXiv:0710.0856 [math.PR]

\bibitem{numericalpde} 
P. Knabner, L. Angermann, {\it Numerical methods for elliptic and 
parabolic partial differential equations}, Texts in Applied Mathematics, 
{\bf 44}, Springer (New York) 2003. 

\end{thebibliography}
\end{document}